\documentclass{amsart}

\usepackage[latin1]{inputenc}
\usepackage{latexsym}
\usepackage{amssymb}
\usepackage{amscd}

\newtheorem{fed}{Definition}[section]
\newtheorem{teo}[fed]{Theorem}
\newtheorem*{teo*}{Theorem}
\newtheorem{lem}[fed]{Lemma}
\newtheorem{cor}[fed]{Corollary}
\newtheorem{pro}[fed]{Proposition}

\theoremstyle{definition}
\newtheorem{rem}[fed]{Remark }

\newtheorem{exa}[fed]{Example}

\newtheorem{num}[fed]{}

\def\fsc{F (S, \c )}
\def\fsa{F (S, \a )}
\def\fcom {the pair $(S, \c )$ is frame admissible}
\def\sii{ if and only if }
\def\linf{\ell^\infty(\N )}
\def\luno{\ell^1(\N )}
\def\lunop{\ell^1(\N )^+}
\def\lunor{\ell_{_\R}^1(\N )}
\def\linfp{\ell^\infty(\N ) ^+}
\def\linfr{\ell_{_\R}^\infty(\N )}
\def\QED{\hfill { $\square$}  }
\def\EOE{\hfill { $\triangle$}  }
\def\eps{\varepsilon}

\def\la{\lambda}
\def\co{\mathrm{conv}}

\def\cF{\mathcal{F}}
\def\H{\mathcal{H}}
\def\N{\mathbb{N}}
\def\R{\mathbb{R}}
\def\LH{\op}
\def\M{\mathbb{M}}
\def\F{\mathcal{F}}

\def\uhs{ \ou{\H}{S} }
\def\uks{ \ou{\cK}{S_1} }
\def\a{\mathbf a}
\def\c{\mathbf c}

\def\s{\mathbf s}
\def\b{\mathbf b}
\def\cA{\mathcal{A}}
\def\cB{\mathcal{B}}
\def\cC {\mathcal{C}}
\def\C{\mathbb{C}}
\def\api{\langle}
\def\cpi{\rangle}

\def\cE{\mathcal{E}}
\def\cK{\mathcal{K}}
\def\cP{\mathcal{P}}

\def\cS{\mathcal{S}}
\def\cM{\mathcal{M}}
\def\cU{\mathcal{U}}
\def\cX{\mathcal{X}}
\def\ese{\cS}
\def\hil{\H}
\def\bm{\left(\begin{array}}
\def\em{\end{array}\right)}
\def\ben{\begin{enumerate}}
\def\een{\end{enumerate}}
\def\beq{\begin{equation}}
\def\eeq{\end{equation}}
\def\barr{\begin{array}}
\def\earr{\end{array}}
\def\bdem{\begin{proof}}
\def\edem{\renewcommand{\qed}{\hfill {{ $\square$}  }}\end{proof}}
\def\calk{\cA (\H )}
\def\inv{^{-1}}
\def\uh{{\mathcal U}(\H )}
\def\uk{{\mathcal U}(\K )}
\def\rai{^{1/2}}
\def\orto{^\perp}

\def\inc{\subseteq}

\def\glh{\mathcal{G}\textit{l} \, (\H)}

\def\ca{L(\H ) }

\def\cak{L(\K ) }

\def\K{{\mathcal K}}

\DeclareMathOperator{\tr}{tr}

\DeclareMathOperator{\clau}{cl}

\newcommand{\posop}{L(\mathcal{H})^+}
\newcommand{\kop}{L_0(\mathcal{H})}

\newcommand{\pint}[1]{\displaystyle \left \langle #1 \right\rangle}
\newcommand{\gen}[1]{\mbox{span}\left\{#1\right\}}

\newcommand{\cluno}[1]{\clau_{\|\cdot\|_1} \left(#1\right)}
\newcommand{\clinf}[1]{\clau_{\|\cdot\|_{_\infty}} \left(#1\right)}
\newcommand{\cln}[1]{\clau_{_{\|\cdot\|}} \left(#1\right)}
\newcommand{\peso}[1]{ \quad \text{ #1 } \quad }
\newcommand{\ou}[2]{ \cU_{#1} (#2) }
\newcommand{\cene}{\mathbb{C}^n}
\newcommand{\mat}{\mathcal{M}_n (\C) }
\newcommand{\matsa}{\mathcal{M}_n(\C)_{h} }
\newcommand{\matpos}{\mathcal{M}_n(\C) ^+}
\newcommand{\matinv}{\mathcal{G}\textit{l}\,_n(\C) }
\newcommand{\matu}{\mathcal{U}(n)}

\newcommand{\bon}{\{e_k\}_{k\in \mathbb{N}}}
\newcommand{\onb}{orthonormal basis }
\newcommand{\bonn}{\{e_n\}_{n\in \mathbb{N}}}

\newcommand{\bonnM}{\{e_n\}_{n\in \M}}
\newcommand{\op}{L(\mathcal{H})}
\newcommand{\oph}{L(\mathcal{H})_h}
\newcommand{\fram}{\{f_n\}_{n\in \mathbb{N}}}
\newcommand{\framk}{\{f_k\}_{k\in \mathbb{N}}}

\newcommand{\framkm}{\{f_k\}_{k\in \M}}

\newcommand{\conv}{\xrightarrow[n\rightarrow\infty]{}}
\newcommand{\convk}{\xrightarrow[k\rightarrow\infty]{}}
\newcommand{\convm}{\xrightarrow[m\rightarrow\infty]{}}

\newcommand{\nulli}[2]{\rm{Null}(#1 , #2)}

\newcommand{\eleuno}{L^1(\hil)}

\begin{document}

\title[ The Schur-Horn theorem for operators and frames]{ 
The Schur-Horn theorem for operators and
frames with prescribed norms and frame operator.}
\author{J. Antezana}
\address{J. Antezana, M.Ruiz, D.Stojanoff, Departamento de Matem\'atica, Universidad Nacional de La Plata, 50 y 115 (1900), La Plata, Argentina and IAM-CONICET, Saavedra 15 (1083), Buenos Aires, Argentina.}
\email{antezana@mate.unlp.edu.ar}
\email{mruiz@mate.unlp.edu.ar}
\email{demetrio@mate.unlp.edu.ar}

\author{P. Massey}
\address{P.Massey,  Departamento de Matem\'atica, Universidad Nacional de La Plata, 50 y 115 (1900), La Plata, Argentina. }
\email{massey@mate.unlp.edu.ar}

\author{M. Ruiz}

\author{D. Stojanoff}

\thanks{Partially supported by CONICET (PIP 2083/00), UNLP (11 X350)
and ANPCYT (PICT03-9521)}


\begin{abstract}
Let $\mathcal H$ be a Hilbert space. Given a bounded positive definite operator $S$ on
$\mathcal H$, and a bounded sequence $\mathbf{c} = \{c_k \}_{k \in \mathbb N}$
of non negative real numbers,  the pair $(S,  \mathbf{c})$ is frame admissible, if there exists a frame
$\{ f_k \}_{k \in \mathbb{N}} $ on $\mathcal H$
with frame operator $S$, such that $\|f_k \|^2 = c_k$,
$k \in \mathbb {N}$. We  relate the existence of such frames
with the Schur-Horn theorem of majorization, and  give a reformulation of the
extended version of Schur-Horn theorem, due to A. Neumann. We use it to
get necessary conditions (and to generalize known sufficient
conditions) for a pair $(S,  \mathbf{c})$,  to be frame admissible.
\end{abstract}

\maketitle

\section{Introduction}
Let $\H$ be a separable Hilbert space and let $S $ be a
bounded selfadjoint operator on $\H$. In the first
part of this note, we give a complete
characterization of the closure in $\linf$ of the set
of possible ``diagonals" of $S$,
i.e., the set $\cC [\uhs ]$ of real sequences  $\c = (c_n)_{n \in \N}$
such that
\beq\label{1}
 \api S e_n , e_n \cpi = c_n \ , \quad  n \in \N\ ,
\end{equation}
for some  \onb $\cB = \bonn$ of $\H$.
Note that, if $\dim \H = m < \infty$, this can be made in terms of
majorization theory. More precisely, the Schur-Horn theorem assures that
$\c \in \R^m$ satisfies Eq. (\ref{1}) for some \onb \sii $\c$ is majorized
by the vector of eigenvalues of $S$ (see Theorem \ref{teo s-h} for a detailed
formulation).
In the general case, we define an analogous form of ``the sum of the greatest
$k$ eigenvalues" in the following way: given $S$ a selfadjoint operator on $\hil$ and $k \in \N$,
we denote
$$U_k(S) = \sup \{ \tr SP : P \in \op
\text{ is an orthogonal projection with } \tr P = k \} \ ,
$$
and $L_k(S) = -U_k(-S)$. We prove, based on the results obtained by A. Neumann in
\cite{neu},  that $\c $ belongs to the $\linf$ - closure of $\cC [\uhs ]$  \sii
\beq\label{infi1}
U_k (\c) \le U_k (S) \peso
{and} L_k(S) \le L_k(\c) \ , \quad k \in \N \  ,
\end{equation}
where $\displaystyle U_k (\c)  =\sup_{|F| = k } \ \sum_{i\in F} c_i \, $, and
$\displaystyle L_k (\c ) =\inf_{|F| = k } \ \sum_{i\in F} c_i = -U_k (-\c )$.
Similarly, if $S$ is a trace class operator, we show that $\c $ belongs to the $\luno$ -
closure of $\cC [\uhs ]$  \sii $\c$ satisfies formulas \ (\ref{infi1}) and $
\sum_{n \in \N} c_n = \tr S$. On the other hand, a somewhat technical characterization
of the maps $U_k$ and $L_k$ is obtained (see Proposition \ref{p:los usubka y eleka de ese}),
which is used  to compute these quantities and to prove their basic properties.
Related results can be found in R. Kadison \cite{Kad1}, \cite{Kad2}, and 
Arveson and Kadison \cite{[AK]} 
(which appeared during the revision process of this 
work).  

In the second part of this note, these extended Schur-Horn theorems are used to give
conditions for the existence of frames with
prescribed norms and frame operator. First we recall some basic definitions.
Let $\M = \N$ or $\M = \{ 1, 2, \dots, m \} := \mathbb {I}_m$, for some $m \in \N$.
A sequence $\framkm$ in $\H$ is called a $frame$
for $\H$ if there exist constants $A,\,B >0$
such that
$$
A\|x\|^2\leq \sum_{k\in \M} |\langle x,f_k \rangle|^2\leq  B\|x\|^2
\ , \peso { for every } x\in \H \ .
$$
For complete descriptions of frame theory and its applications, the
reader is referred to \cite{Dau2}, \cite{[HeiWal]}, \cite{[Ho]},
\cite{[BCHL]} or the books by Young \cite{[Y]} and Christensen
\cite{[liChr]}.
Let $\cF=\framkm $, be a frame for $\H$.
The operator
 \begin{equation}\label{el S}
 S : \H\rightarrow \H \ , \ \   \text{ given by }
 \ \ S(x)=\sum_{k\in \mathbb{M}} \langle x,f_k \rangle f_k  \ ,
 \quad x \in \H \ .
\end{equation}
 is called the \textit{frame operator} of $\F$. It is always bounded,  positive
 and invertible (we use the notation  $S \in \glh^+$).

In the recent works of Casazza and Leon \cite{[cas2]} and
\cite{cas3}, Casazza, Fickus,  Leon and Tremain \cite{[cas1]}, 
Dykema, Freeman, Korleson, Larson, Ordower  and Weber
\cite{Dyk}, Kornelson and Larson \cite{lar2}, and Tropp,
Dhillon, Heath Jr. and Strohmer \cite{[TDHS]},
the problem of existence and (algorithmic)
construction of frames with prescribed norms and frame operator
has been considered. 
Following  \cite{[cas2]}, \cite{cas3}, we  say  the pair $(S, \c )
\in \glh^+ \times \ell^\infty (\M )^+$ is $frame$ $admissible$
if there exists a frame $\cF = \framkm$ on $\H$ such that
\ben
\item $\F$ has frame operator $S$,  and
\item $\|f_k \|^2 =c_k $ for every $ k\in \M$.
\een
In this case, we say that $\F$ is a $(S, \c)-$frame.
We denote  by $\fsc$ the set of all $(S, \c)-$frames on $\H$. Hence
\fcom \ if $\fsc \neq \emptyset \,$.

It is known (see \cite{[cas2]},  \cite{[TDHS]})
that, in the finite dimensional case, there is a connection
between frame admissibility  and the theory of majorization, in
particular with the Schur-Horn theorem. We make this connection
explicit both in the finite and infinite dimensional context. We
use the classical Schur-Horn theorem in the finite dimensional case and its
extension, developed in the first part of the paper, for
the infinite dimensional case.

This presentation of the problem
allows us to get equivalent conditions for the frame admissibility
of a pair $(S, \c) \in \matinv^+ \times \linfp \, $; and necessary
conditions for the frame admissibility
of pairs $(S, \c) \in \glh^+\times \linfp\, $.
We show that, if \fcom , then
$ \sum_{k\in \N} c_k=\infty$,   and
$U_k (\c) \le U_k (S)$ for every $ k \in \N $.
In particular, $\limsup\c  \le \|S\|_{e} $,
the essential norm of $S$ (see Theorem \ref{cor neu frames}).
Then, by strengthening  these conditions we get sufficient
conditions for the frame admissibility of pairs $(S, \c) \in
\glh^+\times \linfp\, $  (Theorem \ref{el teo}). These conditions
are less restrictive that those found by Kornelson and Larson  in
\cite{lar2}.

We briefly describe the contents of the paper. In section 2
we fix our notation, and we state the classical Schur-Horn
theorem. In section 3 we prove the extension of the Schur-Horn theorem
for general selfadjoint operators. In section 4 we give some
reformulations of the notion of frame admissibility which allows us
to apply majorization theory to this problem, and we show equivalent
conditions for frame admissibility in the finite dimensional case
(both for finite or infinite sequences $\c$). In section 5 we
study the infinite dimensional case, showing separately necessary
and sufficient conditions for frame admissibility. In section 6 we
give several examples for the boundary cases of the conditions
studied before. These examples show that, in general, the conditions
can not be relaxed further. We also study different types of frames in
$\fsc$, in terms of their excesses.

\section{Notations and preliminaries.}
Let $\hil$ be a separable Hilbert space, and
$\ca$ be the algebra of bounded linear operators on $\hil$.
We denote $\kop$ the ideal of compact operators,
$\glh$ the group of invertible operators,
$\oph $ the set of hermitian operators,
$\posop$ the set of non negative definite operators,
$\uh$ the group of unitary operators,  and
$\glh^+$ the set of invertible positive definite operators.
We denote by  $\eleuno$ the  ideal of trace class operators
in $\op$. We denote $\eleuno_h = \eleuno\cap \oph$
and $\eleuno^+ = \eleuno \cap \posop$.
We denote by $\luno$ the Banach space of complex absolutely
summable  sequences. By $\lunor$ (resp. $\lunop$)
we  denote the subsets of real (resp. non negative) sequences.
Similarly, we use the notations $\linf$, $\linfr$ and $\linfp$ for bounded sequences.

Given an operator $A\in \ca$, $R(A)$ denotes the range of $A$, $\ker A$ the
nullspace of $A$, $\sigma (A)$ the spectrum of $A$, $A^*$ the adjoint
of $A$, $\rho(A)$ the spectral radius of $A$,  and $\|A\|$ the
spectral norm of $A$.
We say that $A$ is an isometry (resp. coisometry) if $A^*A = I$ (resp. $AA^*=I$).

 We also consider the quotient $\calk = \op/\kop$, which is a
unital $C^*$-algebra, known as the Calkin algebra. Given $T\in \op$,
the {\it essential spectrum} of $T$, denoted by $\sigma_{e}(T)$,
is  the spectrum of the class $T + \kop $ in the algebra $\calk $.
The essential norm $\|T\| _{e} = \inf \{ \|T+K \| : K \in \kop\}$
of $T$ is the (quotient) norm of $T+\kop$, also in $\calk$.
Given $S\in\oph\, $, we define \beq\label{alfas} \alpha^+(S)=\max \sigma_{e}(S) =
\|S\|_{e} \peso{ and } \alpha_-(S)=\min \sigma_{e}(S) \ .
\end{equation}
If $S=\int_{\sigma(S)} t\ dE(t)$ is the
spectral representation of $S$ with respect to the spectral measure
$E$, we shall often consider the following compact operators:
\begin{align}\label{particion del S}
&S^+=\int_{[\alpha^+(S),\ \|S\|]\,}(t-\alpha^+(S))  dE(t) \ , \peso{and} \nonumber \\
&S_-=\int_{[-\|S\|,\alpha_-(S)]}(t-\alpha_-(S))  dE(t) \ .
\end{align}
Note that $S_-\le 0 \le S^+$.

Given a subset $\cM$ of a Banach space $(\cX, \| \cdot \| )$,
its closure is denoted by $\overline{\cM}$ or $\cln {\cM }$,
and the convex hull of $\cM$ is denoted by
$conv(\cM)$. Also, given a closed subspace $ \ese $
of $ \hil$, we denote by $P_\ese$ the orthogonal (i.e.
selfadjoint) projection onto $\ese$. If $B \in \ca $ satisfies
$P_\ese B P_\ese = B$, in some cases we shall use the compression of $B$ to
$\ese$, (i.e. the restriction of $B$ to $\ese$ as a linear
transformation from $\ese$ to $\ese$), and we say that {\it we
consider $B$ as $acting$ on $\ese$}.

Finally, when $\dim \H = n <\infty$,  we shall identify $\H $ with $\cene$,
$\op $ with $\mat$, and we use the following notations:
$\matsa $ for $\oph$, $\matpos$ for $\posop$, $\matu$ for $\uh$,  and $\matinv$ for $\glh$.

\subsection* {Majorization.}
In this subsection we present some basic aspects of majorization theory.
For a more detailed treatment of this notion see \cite{horn2}.
Given $\b=(b_1,\ldots,b_n)\in \R^n$, denote by  $\b^ \downarrow \in \R^n $ the
vector obtained by rearranging the coordinates of $\b$ in non increasing order.
If $\b,\, \c\in \R^n$ then we say that $\c$ is majorized by $\b$,
and write $\c\prec  \b$,  if
$$
\sum_{i=1}^k b^\downarrow_i\geq \sum_{i=1}^k c^\downarrow_i \
\ \ \ k = 1,\dots,n-1, \quad \text{ and } \quad \sum_{i=1}^n b_i=
\sum_{i=1}^n c_i \ .
$$
Majorization is a preorder
relation in $\R^n$ that occurs naturally in matrix analysis.

\begin{fed}\rm
Let $\M = \N$ or $\M = \{ 1, 2, \dots, m \} := \mathbb {I}_m$, for some $m \in \N$.
Let $\cK$ be a Hilbert space with $\dim \cK =
|\M |$ and let  $\cB = \bonnM$  be an orthonormal basis of $\cK$.
\ben
\item For any  $\a = (a_n )_{n \in \M } \in \ell^\infty (\M )$,
denote by $M_{\cB, \a}\in \cak $ the diagonal operator given by
$M_{\cB, \a} e_n = a_n e_n$, $n \in \M$. When it is clear which
basis we are using, we abbreviate $M_{\cB, \a} = M_{\a}$.
\item In particular, for $\a \in \cene$, we denote by $M_{\a} \in \mat $
the diagonal matrix (with respect to the canonical basis of $\cene$)
which has the entries of $\a$ on its diagonal.
 \item The diagonal pinching $\cC _ \cB  :\cak \rightarrow \cak $ associated
 to the basis $\cB$, is defined  by
 $\cC _ \cB (T) = M_{\cB, \a}$, where $\a = (\api Te_n , e_n \cpi )_{n \in \M}$. \EOE
  \een
\end{fed}

\begin{teo}[Schur-Horn] \label{teo s-h}
Let $\b,\ \c\in
\R^n$. Then  $\c \prec  \b $ if and only if  there exists
$U\in \matu$ such that
$$
\cC _ \cE (U^* M_{\b}U)=M_{\c} \ ,
$$
where $\cE$ is the canonical basis  of $\cene$. \QED
\end{teo}

\section{Schur-Horn theorem for selfadjoint operators.}

In this section  we  present a different version of the ``infinite
dimensional Schur-Horn theorem" given by A. Neumann in \cite{neu}. Our
approach avoids the somewhat technical distinction between the
diagonalizable and non diagonalizable case. On the other hand,
this version can be applied more easily to the problem of frame
admissibility in  the infinite dimensional case. The main tools we use
are the Weyl von Neumann theorem and the known properties of approximately
unitarily equivalent operators.

Given a sequence $\a \in \linfr$, Neumann \cite{neu} defines:
\begin{align*}
{U}_k(\a)=\sup_{|F| = k} \ \sum_{i\in F} a_i &&\mbox{and}&&
{L}_k(\a)=\inf_{|F| = k} \ \sum_{i\in F} a_i.
\end{align*}
This generalizes the partial sums which appear in the definition
of majorization. In the first part of this section we shall
extend this definition for arbitrary selfadjoint
operators on a Hilbert space $\hil$.
Denote by $\cP_k$ the set of orthogonal
projections onto $k$-dimensional subspaces of $\H$.

\begin{fed}\label{nuestra def}\rm
Given $S\in \oph$, we define, for any $k \in \N$,
$$
U_k(S)=\sup_{P\in \cP_k}\tr (SP) \quad \text{ and }  \quad
L_k(S)=\inf_{P\in \cP_k}\tr (SP)= - \, U_k(- \,S) \ .
$$
\EOE
\end{fed}

\begin{rem}\label{props uk y lk}\rm It is easy to see that $U_k$ and $L_k$ satisfy the following properties:
\ben
\item For every $k \in \N$,  $U_k$ is a convex map, and $L_k$ is a concave map.
\item The maps $U_k$ and $L_k$ are unitarily invariant,
for every $k \in \N$,  i.e, $U_k(S)=U_k(U^*SU)$, for every $U\in \cU (\H )$ and
$S \in \oph \ $.  \EOE
\een
\end{rem}

The following result asserts that Definition \ref{nuestra def} extends the
natural extrapolation of Neumann's definition for diagonalizable operators.

\begin{pro}\label{ayS}
Let $\cB=\bonn$ be an orthonormal basis of a Hilbert space $\hil$. If $\a\in\linfr$ then, for every $k\in \N$,
\begin{align*}
{U}_k(M_{\cB,\a})={U}_k(\a ).
\end{align*}
\end{pro}
In order to prove this Proposition we need the following
technical results.
\begin{lem}\label{simon}
Let $S \in \kop^+$, and denote by $\la_1 \ge \la_2 \ge  \dots \ge \la_n \ge \dots $
the positive eigenvalues of $S$, counted with multiplicity
(if $\dim R(S) < \infty$, we complete
this sequence with zeros). Then, for every $k \in \N$,
$$
U_k (S) = \sum_{i=1}^k \la_i \ .
$$
Moreover, if $P\in \cP_k$ is the projection onto the subspace
spanned by an orthonormal set of eigenvectors of $\la_1, \dots, \la_k \, $,
then  $U_k (S) =\tr(SP)$.
\end{lem}
\bdem
Fix $k \in \N$.
It suffices to show that
$\tr (SQ ) \le \tr(SP)= \sum_{i=1}^k \la_i   $ for every $Q \in \cP_k$.
This follows from Schur's theorem (the diagonal is majorized by the
sequence of eigenvalues), which also holds in
this setting (see Ch.1 of Simon's book \cite{[Sim]}).
\edem

In \cite{neu}, Neumann proved the following result (Lemma 2.17): if $\a
\in \linfr$,
\begin{align}\label{gene}
a^+_i = \max \{a_i - \limsup\  \a\, ,\ 0\} &&\mbox{and}&&
a^-_i = \min \{a_i - {\liminf} \  \a\, ,\ 0\}  \ ,
\end{align}
$i \in \N$, then, for every $k \in \N$,
\begin{align}\label{con a}
{U}_k(\a)={U}_k(\a^+)+k \  \limsup \ \a &&\mbox{and}&&
{L}_k(\a)={L}_k(\a^-)+k\ \liminf \ \a  \ .
\end{align}
The next result extends Eq. (\ref{con a}) to selfadjoint operators.
This fact is necessary for the proof of Proposition \ref{ayS},
but it is also a basic tool in order to deal with the maps $U_k$ and $L_k\, $.

\begin{pro}\label{p:los usubka y eleka de ese}
Let $S\in \oph$. Then, for every $k \in \N$,
\begin{enumerate}
\item [1.] $U_k(S)=U_k(S^+)+k \  \alpha ^+ (S)$

\item [2.] $L_k(S)=L_k(S_-)+k\ \alpha _-(S) $
\end{enumerate}
where $\alpha ^+(S)$, $\alpha _-(S)$, $S^+$, $S_-$ are defined in
(\ref{alfas}) and (\ref{particion del S}).
In particular,
\beq \label{prome}
\lim_{k\rightarrow \infty}\frac{U_k( S )}{k}=\alpha ^+(S)=\|S\|_{e}  \peso{ and }
\lim_{k\rightarrow \infty} \frac{ L_k(S) }{k}=\alpha_-(S) \ .
\end{equation}
\end{pro}

\bdem
Denote $\alpha^+ = \alpha^+ (S)$, and
\beq\label{P2}
P_2= P_2 (S) = E[\, \|S\|_{e}, \|S\|\, ] = E[\, \alpha^+, \|S\|\, ] \ ,
\end{equation}
where $E$ is the spectral measure of $S$. Recall that
$$
S^+= \int_{[\alpha^+,\ \|S\|]\,}(t-\alpha^+) \ dE(t) =(S-\alpha ^+ )P_2\ .
$$
Then $S-S^+ = S(I-P_2 ) + \alpha^+ P_2 \leq \alpha ^+ I$. Therefore,
for every  $k\in \N $ and $Q\in \cP_k$,
\beq \tr (SQ)=\tr (S^+Q)+\tr ((S-S^+)Q) \leq U_k(S^+)+k \alpha ^+ \ ,
\end{equation}
which shows that  $U_k(S)\leq  U_k(S^+)+k \alpha ^+$ for every $k\in \N $.

To see the converse inequality, suppose first that $\tr P_2 =+\infty$.
Denote by  $\la_1 \ge \la_2 \ge  \dots \ge \la_n \ge \dots $
the eigenvalues of $S^+$, chosen as in Lemma \ref{simon}.

Let $Q_k  \in \cP_k$  be
the projection onto the subspace spanned by
an orthonormal set of eigenvectors of $\la_1, \dots, \la_k \, $.
Then $Q_k \le P_2$. By Lemma \ref{simon},
\[
\tr(SQ_k)  = \tr(S^+ Q_k) + \tr((S-S^+) Q_k) = \sum_{i=1}^k \la _i +k\alpha ^+ =  U_k(S^+)+k\alpha ^+  \ .
\]
Hence, $U_k(S)=U_k(S^+)+k \alpha ^+ $. Now, assume that $\tr P_2
=r <\infty$. If $k\leq r$, the same argument as before shows that
$U_k(S)=U_k(S^+)+k \alpha ^+ $. So, let $k>r$ and take $\eps >0$.
Since $P_{\eps}=E[\, \alpha ^+ -\eps \, ,\, \alpha^+ \, )$ has
infinite rank (otherwise $\|S\|_{e}\leq \alpha ^+ -\eps$), we can
take ${Q} \leq P_{\eps}$ a projection of rank $k-r$. If $Q_k ={Q}+P_2$, 
\begin{align*}
 U_k(S) \geq \tr (SQ_k)&=\tr(SP_2)+\tr(S Q)\\
&=\tr (S^+)+r\alpha^+  + \tr(SP_\eps Q)\\
&\geq \tr (S^+)+r\alpha^+  + (k-r)(\alpha^+ -\eps )\\&=
U_k(S^+)+k\alpha^+  - \eps (k-r) \ .
 \end{align*}
 Since $\eps$ is arbitrary, $U_k(S)=U_k(S^+)+k \alpha ^+ $.
The formula for $L_k(S)$ follows applying item 1 to $-S$.
Finally, as $S^+  \in \kop^+$, then its eigenvalues converge to zero.
Hence, by Lemma \ref{simon}, we get that
$\displaystyle \lim_{k\rightarrow \infty} \frac{U_k (S^+)}{k}= 0$ and
similarly for $L_k(S_-)$.
Therefore, Eq. (\ref{prome})  becomes clear.
\edem

 \textit{Proof of Proposition \ref{ayS}}.
It follows using Lemma \ref{simon}, Proposition \ref{p:los
usubka y eleka de ese}, Eq. (\ref{con a}) and the following
apparent identities: if $S = M_{\cB, \a}$, then
\begin{enumerate}
\item $\alpha^+(S) = \limsup \ \a$,  and \ $\alpha_-(S) = \liminf \ \a$ .
\item $S^+ = M_{\cB, \a^+} $ and $S_- = M_{\cB, \a^-}\ $,
\end{enumerate}
where $\a^+$ and $\a^-$ are defined as in Eq. \eqref{gene}.
\QED

\begin{fed}\rm
Let $\H$ be a Hilbert space, $S\in \op$ and $\cB$ an \onb of $\hil$. Then,
\ben
\item $\ou{\H}{S} =\{U^*SU \, : \, U\in \uh\}$.
\item $\cC[ \ou{\H}{S} ] = \big\{ \c \in \linf \ : \ M_{\cB,\ \c} \in \cC_{\cB}
 (\ou{\H}{S} ) \big\}$. \EOE
 \een
\end{fed}

\begin{rem}
Given $S \in \op$, the definition of $\cC[ \ou{\H}{S} ]$
does not depend on the orthonormal basis $\cB$. In fact, if $\cB'$
is another \onb of  $\H$, $U\in \uh $  maps $\cB$ onto $\cB'$, and
$\c\in \linfp$ satisfies $M_{\cB,\ \c} = \cC_{\cB} (T) $ for some
$T \in \ou{\H}{S}$, then
$$
M_{\cB',\ \c} = UM_{\cB,\ \c}U^* = U\cC_{\cB}  (   T ) U^* = \cC_{\cB'}  ( UTU^* )
\in \cC_{\cB'} (\ou{\H}{S} ) \ .
$$
Therefore
$\big\{ \c \in \linf \ : \ M_{\cB',\ \c} \in \cC_{\cB'}
 (\ou{\H}{S} ) \big\}
= \cC[ \ou{\H}{S} ]$. \EOE
\end{rem}

Given a diagonal operator $M_{\a}\in \oph$, Neumann showed that,
if $\c \in \linfr$, the following statements are equivalent
(Corollary 2.18 and Theorem 3.13 of \cite{neu}):
\begin{enumerate}
    \item ${\c} \in \overline{\cC   [\ou{\H}{M_{\a}}]}$.
    \item $U_k(\a)\geq U_k(\c)$ and  $L_k(\a)\leq  L_k(\c)$, $k \in \N$.
\end{enumerate}
Now, our objective  is to generalize this equivalence for every operator
 $S \in \oph$  (via a reduction to the diagonalizable case).
 We need first the following result about
 approximately unitarily equivalent operators.
\begin{lem} \label{nuevo}
Let $S, T \in \oph$. Then  $S \in \cln{\ou{H}{T}}$ \sii
\ $$\cln{\ou{H}{S}}=\cln{\ou{H}{T}}.$$
In this case $U_k(S) = U_k(T)$ and $L_k(S) = L_k(T) $ for every $k \in \N$.
\end{lem}

\bdem If $\{V_n\} _{n\in \N}$ is a sequence in $\uh$ such that
$\|V_n TV_n^* - S \| \conv 0$, then
$$
\|V_n^* SV_n - T \|  = \|V_n^* (S -V_n TV_n^*)V_n \| =
\|V_n TV_n^* - S \| \conv 0 \ .
$$
Hence $\cln{\ou{H}{S}}=\cln{\ou{H}{T}}$. By
Remark \ref{props uk y lk},  $U_k(V_n TV_n^* ) = U_k(T)$
and $L_k(V_n TV_n^* ) = L_k(T)$, for $n, \, k \in \N$.
Fix $k \in \N$ and take $P \in \cP_k$. Then
$$
\tr SP = \lim _{n\to \infty} \tr V_n TV_n^* P \le \lim _{n\to \infty} U_k(V_n TV_n^* ) = U_k(T) .
$$
Hence  $U_k(S) \le U_k(T)$. Similarly $L_k(S) \ge L_k(T) $. The reverse inequalities
follow from the fact that $V_n^*SV_n \conv T$.
\edem

\begin{rem}\label{Berg}\rm
Two operators  $S, T \in \oph$ satisfying the conditions of Lemma \ref{nuevo} are
called \it approximately unitarily equivalent. \rm This equivalence  relation
is deeply studied in the theory of operator algebras. For example,
as a consequence of the Weyl von Neuman theorem, it
is  proved in Davidson's book \cite{Dav} (II.4.4)
that $S, T \in \oph$ are approximately unitarily equivalent \sii \
$\sigma_e(S)= \sigma_e(T)$ and $\dim\ker (S - \la I) = \dim\ker (T - \la I)$
for every $\la \notin \sigma_e(S)$. From this fact it can be deduced
(see the proof of II.4.4 in \cite {Dav}) that,
for every $S\in \oph$, there exist a diagonalizable $D \in \oph$ which
is approximately unitarily equivalent to $S$. \EOE
\end{rem}

\begin{teo}\label{condiciones necesarias en general}
Let $S\in \oph$ and $\c \in \linfr$. 
Then the following conditions are equivalent:
\ben
\item
$\c \in\overline{\cC  [\ou{\H}{S}]}$.
\item
$U_k(S)\geq U_k(\c)$ and $L_k(S)\leq L_k(\c)$ for every $k\in \N$.
\een
If this is the case, then
$\max \sigma_{e}(S)\geq \limsup\c $  and
$\min\sigma_{e}(S)\leq \liminf\c \, $.
\end{teo}
\bdem
The diagonalizable case was proved by Neumann as we mentioned before.
Note that, in order to deduce our formulation from Neumann's result,
we need Proposition \ref{ayS}.
If  $S$ is  not diagonalizable, by Remark \ref{Berg}, there
exists a diagonalizable  operator $D\in \cln{ \ou{\H}S}$. By Lemma \ref{nuevo},
$U_k(D) = U_k(S)$ and $L_k(D) = L_k(S) $ for every $k \in \N$, and
$\cln{ \ou{\H}D}=\cln{ \ou{\H}S}$.
This implies that
$$
\clinf{\cC  [\ou{\H}{D}]} = \clinf {\cC  [\ou{\H}{S}]},
$$
because the map $T \mapsto \cC_\cB (T)$ is continuous  for every \onb $\cB$.
Hence, the general case reduces to the diagonalizable case.
 The final remark follows from the fact that
 \beq\label{suc}
 \limsup \c = \lim_{k \to \infty} \frac{U_k (\c )}{k}
 \peso {and}
 \liminf \c = \lim_{k \to \infty} \frac{L_k (\c )}{k}  \ ,
 \end{equation}
 and Eq. (\ref{prome}).
  \edem

A similar result can be stated for hermitian operators in $\eleuno$ and
sequences in $\lunor$. In this case our result is just an slight
generalization, using our maps $U_k$ and $L_k$, of some results due to Neumann.

\begin{fed}\rm
Let $\Pi $ be the set of all bijective maps on $\N$ and,
for any $k \in \N$, denote
$\Pi_k \inc \Pi$  the set of permutations $\sigma$ such that
$\sigma (n) = n $ for every $n >k$.
Given $\a \in \ell^\infty  (\N)$ and $\sigma \in \Pi $, we denote
\ben
\item $\a_\sigma = (a_{\sigma (1)}, a_{\sigma (2)}, ....)$.
\item $\Pi\cdot \a = \{ \a \, _ \sigma ,\ \sigma \in \Pi\} $,
the orbit of $\a$, under the action of $\Pi$.
\item $\co(\Pi\cdot \a)$, the convex hull of the orbit of $\a$.\EOE
\een
\end{fed}

\begin{num}\label{neu1}\rm
If $\b , \a $ are sequences in $\lunor$,
Neumann \cite{neu} proved that the following statements are equivalent:
\begin{enumerate}
\item $\b \in \cluno{\co(\Pi\cdot \a)}$.
\item
$\; \sum_{k=1}^\infty  \, b_k =\sum_{k=1}^\infty \, a_k \ $ and
$U_k(\a)\geq U_k(\b)$,  $L_k(\a)\leq L_k(\b)$, $k\in\N$.
\end{enumerate}
\end{num}

\begin{pro}\label{caso l1}
 Let  $S \in \eleuno_h$, and $\b \in \lunor$. Then, the following statements are equivalent,
\ben \item \rm $\b \in\cluno{\cC  [\ou{\H}{S}]}$.\it
 \item $U_k(S)\geq U_k(\b),\; $ $L_k(S)\leq L_k(\b)$ for every $k\in \N$, and
 $\displaystyle \sum_{k=1}^\infty  \, b_k =
\tr S$ .
\een
\end{pro}
\bdem $1 \rightarrow 2$. Note that $\cluno{\cC  [\ou{\H}{S}]} \inc \clinf{\cC
[\ou{\H}{S}]}$. Hence, by Proposition \ref{condiciones necesarias en
general},  $U_k(S)\geq U_k(\b) $ and  $L_k(S)\leq L_k(\b)$
for every $k\in \N$. The equality $\sum_{k=1}^\infty  \, b_k =
\tr S$ clearly holds if $\b \in  \cC  [\ou{\H}{S}]$. The general case follows from
the $\luno$ - continuity of the map $\b \mapsto
\sum_{k=1}^\infty  \, b_k\, $.

$2\rightarrow 1$. Let $\a \in \lunor$ and $\cB = \bon$ an orthonormal
basis of $\hil$ such that $S=M_{\cB ,\a}$. By \ref{neu1} and
Prop. \ref{ayS}, it suffices to show that
$\cluno{\co(\Pi\cdot \a)} \subseteq \cluno{\cC [\ou{\H}{S}]}$.

{\bf Claim:} $\, \cluno{\co(\Pi\cdot \a)}=\cluno{\co(\Pi_0
\cdot \a)}$, where  $\Pi_0=\bigcup_{k\in \N} \Pi_k$.
Indeed, it is sufficient to prove that
$\Pi\cdot \a \subseteq \cluno{\co(\Pi_0 \cdot \a)}$.
Given $\sigma \in \Pi$, $\a_{\sigma} \in \Pi \cdot \a$, and $\eps >0$,  take
$N\in \N$ such that $\sum_{k>N}|a_k| < \frac{\varepsilon}{2}\, $ and
$N_0 \in \N$ such that $\sigma\inv (\mathbb {I}_N ) \inc \mathbb {I}_{N_0}$.
There exists   $\sigma_0 \in \Pi_{N_0}$   such that
$\sigma(k)=\sigma_0(k)$ for every $k \in \mathbb {I}_{N_0} $ such that
$\sigma (k) \in \mathbb {I}_N $. Therefore,
 \begin{align*}
 \| \a_{\sigma}- \a_{\sigma_0}\|_1 &= \sum_{ \sigma (k) 
 \notin \, \mathbb {I}_N }|a_{\sigma (k)}- a_{\sigma_0
 (k)}|\\
 &\leq  \sum_{ \sigma (k)\notin \, \mathbb {I}_N }|a_{\sigma (k)}|
 +\sum_{ \sigma _0(k)\notin \, \mathbb {I}_N }|a_{\sigma
 _0(k)}|<\varepsilon \ .
\end{align*}
 Consider $\b \in \co(\Pi_0 \cdot \a)$. Then, there exists  $n \in
\N$ such that $\b \in \co(\Pi_{n}\a)$.
This means that the first $n$ entries of $\b$ form a  convex
combination of permutations of the first $n$ entries
 of $\a$, and $b_k=a_k$ for every $k> n \, $.
 Hence $(b_1,\ldots , b_{n})\prec(a_1, \ldots , a_{n} )$.
Denote $\cB_n = \{e_k :  k \leq n \}$ and  $\hil_{n} = \gen {\cB_n }$.
So, by Schur-Horn Theorem \ref{teo s-h}, there exists a  unitary  $U_0 \in
L(\hil_n )$ such that
$$
M_{\cB ,\b}|_{\hil_n}= \cC_{\cB_n}(U_0^* M_{\cB, \a}|_{\hil_n} U_0 ) \ .
$$
Letting
$U=\bm{cc} U_0 & 0\\0& I\em  \barr{c} \hil_n \\\hil_n ^\perp \earr \in \uh $,
 we get that
$M_{\cB ,\b}= \cC_\cB(U^* M_{\cB, \a} U )$, and $ \b\in \cC
[\ou{\H}{S}] $. Therefore
$$
\cluno{\co(\Pi\cdot \a)}= \cluno{\co(\Pi_0 \cdot \a)}\subseteq
\cluno{\cC [\ou{\H}{S}]},
$$
which completes the proof.
 \edem

\begin{rem}\rm \label{convexa}
Comparing \ref{neu1} with Proposition \ref{caso l1}, it follows that,
if $S = M_{\cB , \a} $ for some $\a \in \lunor$ and some \onb $\cB$ of $\H$, then
$$\cluno{\co(\Pi\cdot \a)} =  \cluno{\cC [\ou{\H}{S}]}.
$$
In particular,
$\cluno{\cC [\ou{\H}{S}]}$ is a convex set.
On the other hand, since the maps $U_k$ are convex and the maps
$L_k$ are concave, $k \in \N$, it can be deduced from Theorem
\ref{condiciones necesarias en general} that   $\clinf{\cC [\ou{\H}{S}]}$
is convex, for every  $S \in \oph$.
Actually, this fact is known, and can also be deduced from the following
results of Neumann \cite{neu}:
\ben
\item [1.]  If $S = M_{\cB , \a} $ for some $\a \in \linfr$ and
some \onb $\cB$, then $\clinf{\co(\Pi\cdot \a)} =  \clinf{\cC [\ou{\H}{S}]}$.
\item [2.] If $S$ is not diagonalizable, 
\begin{align}\label{tres}
\overline{\cC   [\ou{\H}{S}]} =\overline{\cC   [\ou{\H}{S^+}]}
+[\alpha_-(S),\alpha^+(S)]^\N+\overline{\cC   [\ou{\H}{S_-}]},
\end{align}
where $\alpha ^+(S)$, $\alpha _-(S)$, $S^+$, $S_-$ are defined in
(\ref{alfas}) and (\ref{particion del S}).
\EOE
\een
Note that formula (\ref{tres}), which holds also for
diagonalizable operators, gives another complete
characterization of $\overline{\cC   [\ou{\H}{S}]}$. It can be used
to give an alternative proof of
Theorem \ref{condiciones necesarias en general},
but it can also be
deduced from the statement of this Theorem,
and Proposition \ref{p:los usubka y eleka de ese}.
\end{rem}

\section{Frames with prescribed norms and frame operator.}

\subsection*{Preliminaries on frames.}

We introduce some basic facts about frames in Hilbert spaces. For
a complete description of frame theory and its applications, the reader
is referred to
Daubechies, Grossmann and Meyer \cite {Dau2}, Aldroubi
\cite {Ald}, the review by Heil and Walnut \cite{[HeiWal]}
or the books by Young \cite{[Y]} and Christensen
\cite{[liChr]}.
\begin{fed} \rm
Let $\F = \fram$ a sequence in a Hilbert space $\H$.
\ben
\item [1.]
$\F$ is called a \textsl{frame} if
there exist numbers $ A,B>0$ such that
\beq\label{frame}
A\|f\|^2\leq \sum_{n\in \mathbb{N}}
|\pint{f , f_n}|^2 \leq B\|f\|^2 \ , \peso{  for every } f \in \H \ .
\end{equation}
\item [2.] The optimal constants  $A, B$ for Eq. (\ref{frame})  are called
the \textsl{frame bounds} for $\F$.
The frame $\F$ is called $tight$ if $A=B$, and \it Parseval \rm  if
$A=B=1$. Parseval frames are also called normalized tight frames.\EOE
\een
\end{fed}

\begin{fed}\rm\label{preframe}
Let $\F= \fram$ be a frame in $\hil$. Let $\cK$  be a
separable Hilbert space. Let $\cB = \{ \varphi_n : n \in \N \}$ be an orthonormal
basis  of $\cK$. From Eq. (\ref{frame}), it follows that
there exists a unique $T \in L(\cK , \H )$ such that
$$
T ( \varphi_n ) = f_n  \ , \quad n \in \N \ .
$$
We shall say that the triple $(T, \cK , \cB)$
is a \textsl{synthesis  (or preframe)
operator} for $\F$. Another consequence of Eq. (\ref{frame})
is that  $T$ is surjective.
\EOE
\end{fed}

\begin{rem} \label{cosas} \rm
Let $\F= \fram$ be a frame in $\hil$ and $(T, \cK , \cB)$ a synthesis operator
for $\F$, with $\cB =  \{ \varphi_n : n \in \N \}$.
\ben
\item [1.] The adjoint $T^*\in  L(\H , \cK )$ of
$T$, is   given by $\displaystyle
T^*(x  ) = \sum_{n\in\N} \api x  , f_n \cpi \varphi_n$,
$x  \in \H$. It is called an \textsl{analysis operator} for $\F$.

\item [2.] By the previous remarks, the operator $S = TT^* \in \posop$,
called the \textsl{frame operator} of $\F$,  satisfies
\beq\label{frop}
S f=\sum_{n\in \N} \pint{f , f_n}f_n  \ ,  \peso{for every}  f\in \H \ .
\end{equation}
It follows from (\ref{frame}) that  $AI \leq S \leq B I\, $.
So that $S \in \glh^+$.
Note that, by  formula (\ref{frop}), the frame operator of $\cF$ does not
depend on the chosen synthesis operator.
\EOE
\een
\end{rem}
\begin{fed}\rm\label{intrinsecal}
Let $\F= \fram$ be a frame in $\hil$. The cardinal number
$$
e(\F ) = \dim  \Big\{ (c_n ) _{n \in \N} \in \ell^2(\N ) : \sum_{n\in \N} \  c_nf_n=0\Big\},
$$
is called the \textsl{excess} of the frame.
Holub \cite{[Ho]} and Balan, Casazza, Heil and Landau \cite{[BCHL]} proved that
$$e(\F) = \sup \{ \ |I| :  I \inc \N \ \hbox { and } \
\{ f_n\}_{n\notin I} \ \hbox{ is still a frame on  } \H \} .
 $$
This characterization justifies the name ``excess of $\F$".
It is easy to see that, for every synthesis operator $(T, \cK, \cB )$ of $\F$,
$e(\F ) = \dim \ker T$. The frame $\F$ is called a \textsl{Riesz basis} if $e(\F ) = 0$, i.e., if
the synthesis operators of $\F$ are invertible.
\EOE
\end{fed}
\subsection*{Reformulation of frame admissibility}
 Recall that, given a sequence
 $\c = (c_k)_{k\in \M}\in \ell^\infty (\M )^+$  and $S\in \glh^+ $, we denote by
$\fsc$ the set of $(S, \c)$-frames, i.e., those frames $\cF = \{f_k\}_{k\in \mathbb{M}}$
for $\H$,  with frame operator $ S $, such that $\|f_k\|^2 =c_k$, for every
$k\in \M$, and we say that \fcom \ if $\fsc \neq \emptyset$.
 We shall consider the following equivalent formulation of
frame admissibility, which  makes clear its relationship with the
Schur-Horn theorem of majorization theory.

\begin{pro}\label{la equi}
Let $\c \in \ell^\infty (\M )^+$ and
let $S\in \glh^+$.
Then the following conditions are equivalent:
\ben
\item [1.] The pair $(S , \c ) $ is frame admissible.
\item [2.] There exists a sequence of unit vectors $\{y_k\}_{k\in \M}$
in $\H$ such that
$$
S=\sum_{k\in\M} c_k \, y_k\otimes y_k \ ,
$$
where, if $\M =\N$, the sum converges in the strong operator topology.
\item [3.] There exists
 an extension $\cK = \H\oplus\H_d$ of $\H$
such that, if we denote
\begin{equation}\label{s-h}
S_1 = \begin{pmatrix} S&0\\0&0\end{pmatrix} \barr {c} \H\\\H_d \earr \in
L(\cK ) ^+ \ , \peso{then } \
\c \in \cC  \left[ \ou{\cK}{ S_1}  \right] \ .
\end{equation}
 \een
In this case, there exists a frame $\cF \in \fsc $ with $e(\cF ) = \dim \H_d\ $.
\end{pro}
\bdem The equivalence between conditions 1 and 2 is well known
(see, for example, \cite{Dyk}). Hence we shall prove 1 $
\leftrightarrow$ 3. Assume that $\cF = \{f_k\}_{k\in \mathbb{M}}
\in \fsc \ $. Let $(T_0, \cK_0 , \cB_0 )$ be a synthesis operator
for $\cF$. Consider the  polar decomposition $T_0=U|T_0|$, where
$U:\cK_0 \rightarrow \H$ is a coisometry with initial space
$(\ker T_0)^\perp$ and range $\H $. Note that $U^*$ maps
isometrically $\H$  onto $\ker T_0 \orto $. Denote $\H_d = \ker
T_0$, and $\cK = \H \oplus \H_d \ $. Let $V : \cK \to \cK_0$ be
the unitary operator given by

$$
 V (\xi _1 , \xi_2 ) =   U^* \xi_1
 + \xi_2  \ \  , \quad \peso {for }
(\xi _1 , \xi_2 ) \in \H \oplus \H_d =   \cK \ .
$$
Consider the \onb  $\cB = V^*(\cB_0)$ of $\cK$, and
$T = T_0V \in L(\cK , \H )$. Then $(T, \cK , \cB )$ is another
synthesis operator for $\cF$, with $\ker T = \H_d$.

Let  $ T_1 \in L(\cK )$
given by $ T_1\xi=T\xi\oplus 0_{\H_d}$, \ $\xi \in \cK $. Then
$ T_1 ^* T_1=T^* T$, $|T_1| = |T|$, and
$$
T_1 T_1 ^*= \begin{pmatrix} T T^* & 0\\ 0 &0\end{pmatrix} \barr{c} \H \\\H_d \earr
=\begin{pmatrix} S& 0\\ 0 &0\end{pmatrix}   = S_1 \ .
$$
If $T_1 = U_1 |T_1|= U_1 |T|$ is the
polar decomposition of $T_1$, then $U_1$ acts on $\H = (\ker T_1)\orto$
as a unitary operator.
Hence $W = U_1 + P_{\H_d} \in \uk $.
Since  $T_1= W |T|$,
$$
S_1 =  T_1 T_1^* = W  |T|^2  W  ^*= W (T^*T)W ^* \quad \implies \quad
  W ^* S_1 W =T^*T \ .
$$
On the other hand, if $\cB = \bon$, then
$\langle T^*T e_k,e_k\rangle= \langle T e_k,Te_k\rangle=\|f_k\|^2=c_k$,
for every $k \in \M$.
Therefore,
$$
\cC _ \cB  \left(  W ^* S_1  W  \right) =
\cC _ \cB  (T^*T) = M_{\cB , \c} \ \  \implies  \ \
\c \in \cC  \left[ \ou{\cK}{ S_1}  \right]\ .
$$
Conversely, suppose that
there exists an extension $\cK = \H\oplus\H_d$ of $\H$ and
$ V \in \uk $ such that $M_{\cB , \c} = \cC_\cB (V^*S_1 V)$,
for some \onb $\cB = \bon$ of $\cK$.
Let $T = S_1 \rai V$. Since $S$ is invertible, then  $R(T) = \H$ and
$\dim \ker T = \dim \H_d$.
Thus $\cF = \{Te_k\}_{k\in \M}$ is a frame for $\H$, with
frame operator $TT^*\big|_{_\H} = S_1 \big|_{_\H} = S $.
Since $T^*T=  V ^* S_1 V \ $ and $  \cC_\cB (V^*S_1 V) = M_{\cB , \c} \ $,
then $\|Te_k \| ^2 = \api T^*T e_k , e_k \cpi = c_k $, for every $k \in \M$.
Hence $\F \in \fsc$ with $e(\cF ) = \dim \H_d \, $.
 \edem

\subsection*{The finite-dimensional case}
In this section  we assume that $\H$ is finite dimensional. We shall consider separately,
the cases of  frames of finite or infinite length.
Suppose that $S \in \matpos$ and $|\M| = m <\infty$. In this case, the
classical Schur-Horn Theorem \ref{teo s-h} gives  a complete characterization of
frame admissibility for  $(S ,\c)$.

 \begin{teo} \label{cor lar}
 Let $\c\in  \R_{>0}^m$ and let $S\in
\matinv^+$, with eigenvalues
$b_1\geq b_2\geq\ldots\geq b_n>0$. Then, \fcom \
 if and only if $$\sum_{i=1} ^k b_i \geq
\sum_{i=1} ^k c_i \ \ \text{ for }\ \ 1\leq k\leq n-1 \ , \quad
\text{ and } \quad  \sum_{i=1} ^n b_i=\sum_{i=1} ^m c_i \ .
$$
In other words, if \ $\c \prec (b_1 ,  \dots , b_n ,0 ,\dots  , 0 ) \in \R^m$.\QED
 \end{teo}

This result was obtained in \cite{[cas2]} and \cite{lar2}, from an
operator theoretic point of view. Actually the proofs given there
can be adapted so as to obtain a proof of the classical Schur-Horn theorem
that are quite conceptual and simpler than those in the
literature.
Now, we consider frame admissibility  for infinite  sequences
in finite dimensional Hilbert spaces. The case $S = I$ of the next
result appeared in  \cite{[cas1]}. 

\begin{teo}\label{infinitos en finito}
Let $\c \in \ell^\infty (\N )^+$.
Let  $S\in \matinv^+$, with eigenvalues
$b_1\geq b_2\geq\ldots\geq b_n>0$.
  Then the following conditions are equivalent:
  \ben
  \item \fcom .
 \item
 $ \sum_{i=1} ^k b_i \geq U_k(\c) ,\text{ for every }1\leq k\leq
n-1 \text{, and }\sum_{i=1} ^n b_i=\sum_{i\in \N} c_i \, .$
 \een
 \end{teo}

\bdem  Let $\b =(b_1,\ldots, b_n,0,\ldots,0,\ldots) \in \linfp$.

$2 \rightarrow 1$:
Let $\H$ be a  infinite dimensional Hilbert space,
and consider
$$
 S_1 = \begin{pmatrix} S&0\\0&0\end{pmatrix}\in L(\cene \oplus\H ) \ .
$$
Then there exists an orthonormal basis $\cB = \bon$
of $\cK = \cene \oplus\H $
such that $ S_1  = M_{\cB, \b} \, $. Hence, by Proposition \ref{ayS},
$$
U_k( S_1 )= \sum_{i=1}^{k} b_i\geq U_k(\c),\
\ \text{for every } k\in \N.
$$
On the other hand, note that
$L_k( S_1 )=0 \le L_k(\c)$ for every $k\in \N$ and
$\sum_{i=1}^n b_i =\sum_{i\in \N} c_i$.
Then, by  Proposition \ref{caso l1},
there exists a sequence $\{V_m \} _{m\in \N} $ in $\uk$  such that
$$
\cC _ \cB  \left(V_m ^* S_1 V_m\right )\xrightarrow[m\to \infty ]{\|\, \|_1}
M_{\c} \ ,
$$
where $\|A\|_1=\tr |A|$. Therefore, by Proposition \ref{la
equi}, there exists a norm bounded sequence of epimorphisms
$T_m: \cK\rightarrow \cene$ such that that $T_mT_m^*=S$ for all
$m\in \N$,  and $(\|T_m( e_i)\|^2)_{i\in
\N}\xrightarrow[m\to \infty]{\luno}\c$. Then, by a standard
diagonal argument, we can assure the existence of a subsequence, which we still
call $\{T_{m}\}_{m\in \N}$,
such that
$$
T_{m}(e_i)\xrightarrow[m\to \infty]{}f_i\in \cene , \ \text{ with }
\|f_i\|^2= c_i \ \text{ for every } i\in \N.
$$
Let $ T_0 : \gen{\cB} \to \cene $ be the unique (densely defined) operator,
such that $T_0(e_i)=f_i$ for every $i\in \N$.
Note that $T_0$ is bounded because, if
$x=\sum_{i=1}^r \alpha_i \, e_i$ and $C = \sum_{i\in \N} c_i = \tr S$,  then
\begin{eqnarray*}
\|T_0(x)\| & =&  \|\sum_{i=1}^r \alpha_i \,f_i\|\leq
\sum_{i=1}^r |\alpha_i| \|f_i\| \\
& \leq & \left( \sum_{i=1}^r c_i\right)^{1/2}
\left(\sum_{i=1}^r|\alpha_i|^2\right)^{1/2}\leq C\rai \|x\| \ .
\end{eqnarray*}
  The bounded extension of $T_0$ to $\cK$ is denoted  $T$.

{\bf Claim :} $\|T_{m}-T\|\convm 0$.

Indeed, let $\eps>0$ and
$i_0\in\N$ such that $\sum_{i= i_0}^\infty  c_i< \eps$. Then, there
exists $m_1\in \N$ such that
\begin{equation}\label{una ec}
\sum_{i= i_0}^\infty  \|T_m(e_i)\|^2\leq \eps,\,\text{
for every } m\geq m_1 \ .
\end{equation}
This is a consequence of the fact that
$(\|T_m(e_i)\|^2)_{i= i_0}^\infty \xrightarrow[m \to \infty ]{\luno}
(c_i)_{i= i_0}^\infty $. On the other hand,
there exists $m_2 \ge m_1$ such that
\begin{equation}\label{dos ec}
\sum_{i=1}^{i_0-1} \|T_m(e_i)- f_i \|^2\leq \eps,\,
\text{ for every } m\geq m_2\ .
\end{equation}
Let $m \ge m_2$
and $x=\sum_{i=1}^r \alpha_i e_i \in \gen{\cB} $.
By equations (\ref{una ec}) and (\ref{dos ec}),
\begin{align*}
&\|(T_m -\  T)(x)\|^2 \leq
\left(\sum_{i=1}^r| \alpha_i|^2 \right) \left(\sum_{i=1}^r\|
(T_m - T)(e_i)\|^2 \right)\\
& \displaystyle\leq  \|x\|^2
\left( \sum_{i=1} ^{i_0-1}\| (T_m - T)(e_i)\|^2+ 2\sum_{i= i_0}^\infty  \| T_m(e_i)\|^2+
\| T(e_i) \|^2 \right)\\
& \leq  5 \, \eps \, \|x\|^2
\ ,
\end{align*}
which proves the claim. Therefore
$\displaystyle TT^*=\lim_{m\to \infty }T_mT_m^*=S$. We have proved that
the frame $\cF = \{f_i\}_{i\in \N} \in \fsc $.

$1 \rightarrow 2$: This follows from 
Theorem
\ref{condiciones necesarias en general}, applied to $S_1$ and $\c$,
and Proposition \ref{la equi}.
\edem

\begin{rem}
The statement of Theorem \ref{infinitos en finito} can be reformulated in terms
of finite rank operators and sequences in $\luno$ in the following way:
Let $\cK$ be a separable, infinite dimensional Hilbert space. Let $S_1 \in L(\cK )^+$
such that
$\dim R(S_1) < \infty$. Then $\cC[\uks ]$ is closed, as a subset of $\luno$.

Indeed, suppose that $S_1 \neq 0$ (the case $S_1 = 0$  is trivial). Then,
there exists a sequence $\b =(b_1,\ldots, b_m,0,\ldots,0,\ldots) \in \lunop$,
with $b_m >0$, and an \onb $\cB= \bonn$ of $\cK$ such that
$S_1 = M_{\cB, \b}\, $. Let $\c \in \lunop$.
By Proposition \ref{caso l1}, condition 2 of Theorem \ref{infinitos en finito} means that
$\c \in \cluno{\cC[\uks ]}\, $.
But, by Proposition \ref{la equi},
condition 1 of Theorem \ref{infinitos en finito} means that
$\c \in \cC[\uks ]\, $.

Note that, although $\cluno{\co(\Pi\cdot \b)} = \cluno{\cC[\uks ]} = \cC[\uks ]$,
as it is shown in Remark \ref{convexa},
 it is not true that $\co(\Pi\cdot \b)$ is closed, as a subset of $\lunop$.
 For example, if $\b =(1, 0,0,\ldots)$, then, by Proposition \ref{caso l1},
 $$ \c = \big( \frac1{2^{n}}\big)_{n \in \N} \in  \cluno{\cC[\ou{\cK}{e_1\otimes e_1}]}
 = \cluno{\co(\Pi\cdot \b)}\ .
 $$
 Nevertheless,  $\c \notin  \co(\Pi\cdot \b)$,
 because every sequence in $\co(\Pi\cdot \b)$ has finite non zero entries.
In this case, $\c = \cC_\cB ( x \otimes x ) \in \cC[\ou{\cK}{e_1\otimes e_1}]$, where
$ x = \sum _{n \in \N} 2^{-\frac n2} e_n $.
\EOE

\end{rem}

\section{The infinite-dimensional case}

Throughout this section $\H$ denotes a separable infinite
dimensional Hilbert space.
The first result gives necessary conditions for
frame admissibility:

\begin{teo}\label{cor neu frames}
Let  $S\in \glh^+$ and $\c \in \ell^\infty (\N )^+$.
If \fcom , then $\displaystyle  \sum_{i\in \N} c_i=\infty$,   and
$U_k(S)\geq U_k(\c)$, for every $k \in \N$. In particular, $\limsup\c  \le \|S\|_{e} $.
\end{teo}

\bdem
Suppose that there exists  a frame $\cF \in \fsc$.
Then, by Proposition \ref{la equi},
there exists  an extension $\cK = \H\oplus\H_d$ of $\H$
such that, if we denote
$$
S_1 = \begin{pmatrix} S&0\\0&0\end{pmatrix} \barr {c} \H\\\H_d \earr \in
L(\cK ) ^+ \ , \peso{then } \
\c \in \cC  \left[ \uks  \right ] \ .
$$
Hence, $\sum_{i\in \N} c_i = \tr M_{\c} = \tr S_1 = \infty$ .
On the other hand, by Proposition \ref{p:los usubka y eleka de ese},  $U_k (S) =
U_k (S_1)$ for every $k \in \N$. Then,  applying
Theorem \ref{condiciones necesarias en general}, the
statement follows.
\edem

\begin{rem}\rm
Let  $S\in \glh^+$ and $\c \in \ell^\infty (\N )^+$.
Then, by Theorem \ref{condiciones necesarias en general}
and Proposition \ref{la equi},
the following conditions are equivalent :
\ben
\item  $U_k(S)\geq U_k(\c)$ for every $k\in \N$.
\item There exists a sequence
$\mathcal F_k=\{f_{ik} \}_{i\in \N}$, $k\in \N$ of frames on $\H$,
such that $S$ is the frame operator of every $\mathcal F_k$ and
$\|f_{ik}\|\convk  \sqrt{c_i}$ uniformly for $i\in \N$.
\een
Indeed, note that the inequalities involving the maps
$L_k$,
$k\in \N$, can always be fulfilled if we consider a sufficiently
large extension $\H\oplus\H_d$ of $\H$.
In this case, $\limsup\c  \le \|S\|_{e} $. \EOE
\end{rem}

 At this point we should note that the conditions of
Theorem \ref{cor neu frames} are not sufficient  to assure that \fcom , as Example
\ref{contraejemplo} below shows. That is, we can not remove the
closures in the equalities of Theorem \ref{condiciones necesarias en general},
as it was already mentioned in \cite{neu}, for the diagonalizable case.

In \cite{lar2} (see also \cite{[cas1]}) appears the following  result which gives
sufficient conditions for a pair $(S, \c) $
in order to be frame admissible:

\begin{teo}[Kornelson-Larson] \label{teo lar-kor}
Let $S\in \glh^+$ and $\c\in l^\infty(\N)^+$.
Suppose  that $\sum_{i\in \N} c_i=\infty$ and
$\|\c\|_\infty  < \|S\|_{e}$. Then \fcom .
\QED
\end{teo}

The following result, which generalizes Theorem \ref{teo lar-kor},
strengths slightly the necessary  conditions for frame admissibility
given by Theorem \ref{cor neu frames},  to  get sufficient
conditions.  A tight frame version of this result appeared in 
R. Kadison \cite{Kad1} and \cite{Kad2}. 
Recall the notation $P_2(S) = E[\, \|S\|_{e}, \|S\|\, ]$, where $E$
is the spectral measure of $S \in \posop$.

\begin{teo}\label{el teo}
Let $S\in \glh^+ $ and $\c\in l^\infty(\N)^+$,
such that $\sum_{i\in \N} c_i=\infty$.
 Assume one of the following two conditions:
\ben

\item  [1.]
\ben
\item $\tr  P_2(S) =\infty$,
\item $U_k(S)\geq
U_k(\c)$ for every $k\in \N$,  and
\item $ \|S\|_{e}> \limsup(\c).$
\een
\item []
\item [2.]
\ben
\item $\tr  P_2(S) =r\in \N$,
\item $U_k(S)\geq U_k(\c)$ for $1\le k \le r$,
\item $U_{k}(S)>  U_{k}(\c)$, for $k>r $,  and
\item $ \|S\|_{e}>\limsup(\c).$
\een
\een
Then, \fcom .
\end{teo}

\bdem
By Proposition \ref{la equi},  it suffices to show that the there
exists a sequence of unit vectors $\{x_k\}_{k\in \N}$ such that
$
S=\sum_{k\in \N} c_k \ x_k\otimes x_k \ .
$
Assume that the first condition holds.
Then, since $ \|S\|_{e}> \limsup(\c)$, 
there exist $m_0\in \N$ and $\eps>0$ such that
$$
c_m \le\|S\|_{e}-\eps \peso{for $m \ge m_0$ }
 $$
Let $\mu_1\geq \mu_2\ldots\geq \mu_{n} \ge \dots  $ be the sequence of eigenvalues
of $S^+$, chosen as in Lemma \ref{simon}.  Let $\{y_n\}_{n\in\N} $ be an
orthonormal system such that $S^+ y_n = \mu_n y_n$.
Denote $\la_n = \mu _n +  \|S\|_{e}$, $n \in \N$.
Note that $\|S\| \ge \lambda_n \ge \|S\|_{e}$, and  $S y_n = \la_n y_n$,
$n \in \N$.
By Proposition \ref{p:los usubka y eleka de ese}, for every $k \in \N$,
$$
\sum_{i=1}^{k} \lambda_i\, y_i\otimes y_i\leq S \ ,
\peso{ and }
U_k(S)=\sum_{i=1} ^k \lambda_i \ .
$$
Let $n_0$ be the first integer such that
$\displaystyle \sum_{i=1} ^{n_0} c_i>\sum_{i=1}^{m_0} \lambda_i$ .
Then $n_0\geq m_0+1$,  and
$$
h=\sum_{i=1} ^{n_{0}} c_i- \sum_{i=1}^{m_0}\lambda_i \leq
c_{n_0}<\|S\|_{e} \le \la_{m_0+1}\ .
$$
Let $\c_0 = (c_1 , \dots, c_{n_0})$. Since
$$
 \sum_{i=1} ^k \lambda_i=U_k(S)\geq U_k(\c)\geq U_k(\c_0) \ , \quad
1\le k\le m_0 \ ,
$$
then $\c_0 \prec (\la_1 , \dots, \la_{m_0}, h , 0 , \dots , 0) \in \R^{n_0}$.
Denote by
$$
\displaystyle S_1 = h\ y_{m_0+1}\otimes y_{m_0+1}+
\sum_{i=1} ^{m_0} \lambda_i \ y_i\otimes y_i \le S  \ ,
$$
and $S_2 = S-S_1$. Then, the pair
$(S_1, \c_0 \, )$, acting  on  $\gen{y_1, \dots , y_{m_0+1}}$,
satisfies the conditions of Theorem \ref{cor lar}. Hence,
there exists a set of unit vectors $\{x_1,\ldots, x_{n_0}\}$
such that
$\displaystyle
\sum_{i=1} ^{n_0} c_i \ x_i\otimes x_i = S_1
 \  .
$
Note that $S_2 \ge 0$, $R(S_2 ) $ is closed (by Fredholm theory),
and $\|S_2\|_{e}=\|S\|_{e}$.
Then, we can apply Theorem \ref{teo lar-kor} to the pair
$(S_2 , \{c_i\}_{i>n_0})$, acting on $R(S_2 ) $. So,
there exist unit vectors $x_k$, for $k >n_0$, such that
$$
S_2 = \sum_{i=n_0+1}^\infty  c_i \ x_i\otimes x_i \ .
$$
Therefore we obtain the rank-one decomposition $S=\sum_{i\in \N} c_i\ x_i\otimes x_i$.

 Assume condition 2.  
Note that, by equations (\ref{prome}) and (\ref{suc}),  the condition $ \|S\|_{e}> \limsup(\c)$
implies that $U_{m}(S)-U_{m}(\c)\convm \infty$. Therefore, by item (c), we can assume that there exists $\delta>0$ such that
\begin{enumerate}
  \item $U_{r+k}(S)-\delta > U_{r+k}(\c)$, for every $k\in\N$.
  \item There exists $m_0\geq 1$ such that
  $c_m \le\|S\|_{e}-\delta$ for $m \ge m_0$ .

\end{enumerate}
Let $m_1=\max \{ m_0,r+1 \}$.
Let $\mu_1 \ge  \dots \ge  \mu_r$ be the greatest eigenvalues of $S^+$, and let
$\{ y_1 , \dots, y_r \} $ be an associated orthonormal set of eigenvectors.
Denote  by
$$
\la_i = \mu_i + \|S\|_{e} \, , \  1\le i \le r \, , \ \text{ and }
\ \la _i = \|S\|_{e}-\frac{\delta}{2m_1} \, , \ \ r+1 \le i \le m_1+1 \, .
$$
Then, by Proposition \ref{p:los usubka y eleka de ese},
\begin{enumerate}
   \item $\displaystyle U_k(S)=\sum_{i=1} ^k \lambda_i  $ , for  $1\le  k \le r $, and 
  \item $\displaystyle U_k(\c)\leq U_k(S)-\delta\leq
  \sum_{i=1} ^k \lambda_i $ , for $r+1\le k \le m_1+1 $ .
\end{enumerate}
On the other hand, since $Q = E([\|S\|_{e}-\delta/2m_1 , \|S\|_{e}) \, )$
has infinite rank, there exists an orthonormal set
$\{y_{r+1}, \dots , y_{m_1+1}\} \inc R(Q)$. Therefore
$$
\displaystyle \sum_{i=1}^{m_1+1} \lambda_i \ y_i\otimes y_i\leq S \  .
$$

Let $n_0$ be the first
integer such that $\sum_{i=1} ^{n_0} c_i>\sum_{i=1}^{m_1} \lambda_i$.
Then $n_0\geq m_1+1$ and
$$
h=\sum_{i=1} ^{n_{0}} c_i- \sum_{i=1}^{m_1}\lambda_i\leq c_{n_0}\leq\|S\|_{e}-\delta \le \lambda_{m_1+1} \  .
$$
Let $\c_0 = (c_1 , \dots, c_{n_0})$. Since
$$
\ \sum_{i=1} ^k \lambda_i=U_k(S)\geq U_k(\c)\geq U_k(\c_0)
\ ,\quad 1\le k \le r \ , \quad \text{and}
$$
 $$\ \ \ \ \sum_{i=1} ^k \lambda_i\geq U_k(S)-\delta\geq
U_k(\c)\geq U_k(\c_0)   \ ,\quad r+1\le k \le m_1 \ ,
$$
then $\c_0 \prec (\la_1 , \dots, \la_{m_1}, h , 0 , \dots , 0 ) \in \R^{n_0}$.
So, by Corollary \ref{cor lar}, there exists a set of unit vectors
$\{x_1,\ldots, x_{n_0}\} \inc \H$
such that
$$
S_1 = \sum_{i=1} ^{m_1} \lambda_i \
y_i\otimes y_i +h\ y_{m_0+1}\otimes y_{m_0+1}= \sum_{i=1} ^{n_0}
c_i \ x_i\otimes x_i.
$$
Since $S_1
\leq \sum_{i=1}^{m_1+1} \lambda_i \ y_i\otimes y_i $ , then
$
S_2 =S- S_1 \ge 0$
and $\| S_2 \|_{e}=\|S\|_{e}$.
 As before, we apply Theorem
\ref{teo lar-kor} to the pair $(S_2 , \{c_i\}_{i>n_0})$, acting on $R(S_2 ) $,
and we obtain a decomposition
$$
S_2=\sum_{i= n_0+1}^\infty  c_i \ x_i\otimes x_i.
$$
Therefore we obtain
the rank-one decomposition $S=\sum_{i\in \N} c_i\ x_i\otimes x_i$.
\edem

\rm
Example \ref{nada} below shows that the  condition 2 (c)
of Theorem \ref{el teo}
can not be dropped in general.

\begin{cor}
Let $0< A\in \R$ and  $\c\in \ell^\infty (\N )^+$ such that $0<c_i\leq A$,
$i \in \N$. Denote $J=\{i\in \N: c_i=A\}$. Assume that
$$
\sum_{i\notin  J } c_i=\infty \peso{ and }\limsup_{ i\notin  J }  \, c_i < A
\peso{(or, equivalently, } \sup_{ i\notin  J }  \, c_i < A) \ .
$$
Then the pair
$(AI, \c)$ is admissible. This means that there
exists a tight frame with norms prescribed by $\c$ and frame
constant $A$. \QED
\end{cor}

\section{Some examples}
In the following example we shall see that
$$U_k(S) >  U_k(\c) , \ k \in \N \peso{and} \|S\|_{e}=\limsup(\c) \ \ \not\Rightarrow \ \ \fsc \neq \emptyset \ .
$$
\begin{exa} \label{contraejemplo}
Let $S=I\in\LH$ and $a\in(0,1)$. Let $\c \in \ell^\infty (\N )^+$ be given
by $c_1=p\in(0,1)$ and
\[ c_k = \left\lbrace
  \begin{array}{c l}
    a^{k} & \text{if $k\neq 1$ is odd },\\
    1-a^{k}   & \text{if $k$ is even}.
  \end{array}
\right. \] Then, $0<c_k<1$ for $k\in \N$, $\sum_k c_k=\infty  = \sum_k (1-c_k)\ $,
and $\limsup \c = 1 = \|S \|_{e}\ $. Suppose that there exists a frame
$\cF = \framk  \in \fsc $.
Then
$$
\|x\|^2=\sum_{k \in \N} |\langle x,f_k \rangle|^2 \ , \quad  \text{for
every  }x\in\H \ .
$$
In particular, we get, for every $j \in \N$,
$$
\|f_j\|^2=\sum_{k \in \N} |\langle f_j , f_k \rangle|^2=
\|f_j\|^4 + \sum_{k\neq j} |\langle f_j , f_k \rangle|^2.
$$
Thus, if $j\neq 1$, we obtain the inequality
$$
|\langle f_1 , f_j \rangle|^2 = |\langle f_j , f_1 \rangle|^2 \leq \sum_{k \neq j} |\langle
f_j,f_k \rangle|^2=\|f_j\|^2 - \|f_j\|^4=c_j(1-c_j).
$$
Therefore,
\begin{eqnarray}\label{pp}
\nonumber p=\|f_1\|^2 & \leq & \|f_1\|^4 + \sum_{j \neq 1} c_j(1-c_j)
 \\ &=& p^4 + \sum_{j \neq 1} a^j(1-a^j)= p^4 +\sum_{j \neq 1}  a^j - \sum_{j \neq 1} a^{2^j}
 \\ \nonumber &=& p^2 + \frac{1}{1-a} - \frac{1}{1-a^2}= p^2 + \frac{a}{1-a^2}
\end{eqnarray}
Taking $\displaystyle  p=\frac{1}{2}$ and $a\in (0,1)$ such that
$\displaystyle \frac{a}{1-a^2}<\frac{1}{4}$ ,  we get that
$\displaystyle  p> p^2+ \frac{a}{1-a^2}\ $,
contradicting Eq. (\ref{pp}). Hence, in this case,  $\fsc = \emptyset$. Note that
the pair $(S, \c )$ satisfies all necessary conditions of Theorem \ref{cor neu frames},
because $U_k(S) =k =  U_k(\c)$  for every $k \in \N$.
 \EOE
\end{exa}

In the second example we see that, in general,
$$U_k(S) \ge  U_k(\c) , \ k \in \N \peso{and} \|S\|_{e}>\limsup(\c) \ \ \not\Rightarrow \ \ \fsc \neq \emptyset \ .
$$
\begin{exa}\label{nada}
Let $S=M_{\s}$ be the diagonal operator, with respect to an \onb
 of $\H$, given by $\s=\{1-(i+1)\inv\}_{i\in \N}$ and let
$(c_i)_{i\in \N}$ be given by $c_1=1$ and $c_i=1/2$ for every
$i\geq 2$. Note that
\begin{itemize}
\item $1=\|S\|_{e}>1/2=\limsup(\c)$,
\item $U_1(S) = U_1(\c) $, and
\item  $U_k(S)=k > 1 + (k-1)/2 = U_k(\c)$  for every $k \ge 2$.
\end{itemize}
Still, $\fsc = \emptyset$.
Indeed, suppose that there exists $\F \in \fsc$. Then, by Proposition \ref{la equi}
there exists  an extension $\cK = \H\oplus\H_d$ of $\H$
such that, if
$$
S_1 = \begin{pmatrix} S&0\\0&0\end{pmatrix} \barr {c} \H\\\H_d \earr \in
L(\cK ) ^+ \ , \peso{then } \
\c \in \cC  \left[ \uks \right ] \ .
$$
Let $V \in \uk$ such that, in a \onb $\cB = \bon$, $M_\c = \cC_\cB (V^*S_1 V)$.
Take $x= P_\hil Ve_1$. We have that $\|x\|\le 1$ and
$\api Sx, x\cpi = \api M_{\c} e_1, e_1\cpi = c_1 = 1$, while $\|S\| = 1$.
Then $Sx = x$,  and $1$ would be an eigenvalue of $S$, which is false.
In this example, condition 2 (c) of Theorem \ref{el teo}
does not hold,   because
$\|S\|=\|S\|_{e}$, which implies that $r = \tr P_2(S)=0$;
but  $U_1(S)=1=U_1(\c) $ .
Note that 
$\sum_k c_k=\infty  = \sum_k (1-c_k)$, as in the previous example.
 \EOE
\end{exa}

\subsection*{The excess of frames in $\fsc$}

Let  $S\in \glh^+$ and $\c=(c_i)_{i\in \M} \in \ell ^\infty (\M )^+$
such that \fcom .
Then, there can be many different types of
frames $\cF  \in \fsc$. We consider the set
$$
\nulli{S}{\c} =\{ \ e(\cF) : \cF \in \fsc \ \} \ .
$$
In the Example below, we show that this set can be as big as possible.
Moreover, this example shows that there exists an admissible
pair $(S, \a )$, satisfying just the necessary conditions of
Theorem \ref{cor neu frames} and in this case $U_k (S) = U_k (\a )$, $k \in \N$, and
$\limsup \a = \|S\|_{e}\, $.

\begin{exa}\label{grande}
Let $\H$ be a Hilbert space with an orthonormal basis
$\cB = \{x_n\}_{n \in \N}$. Let
$$
\a = \Big(\frac12 \, , 1, \frac12 \, , 1, \frac12 \, , \dots \Big) \in \ell^\infty (\N)^+ \ ,
\peso{ and }  S = M_{\cB , \a} \in \glh^+ \ .$$
Then, the frame (Riesz basis) $\F_0 = \{a_n\rai x_n\}_{n\in \N}$ has frame operator $S$,
so that $\F_0 \in \fsa$.
On the other hand, let
$$
\F_1 = \Big\{ \frac1{\sqrt{2}} \,  x_2 , x_4, \frac1{\sqrt{2}} \, x_2 \, , x_6 ,
 \frac1{\sqrt{2}}\, x_1 , x_8 ,  \frac1{\sqrt{2}}\, x_3 , x_{10} , \dots \Big\} \ .
$$
It is easy to see that also $\F_1 \in \fsa$, but $e(\F_1 ) = 1$. Analogously,
$$
\F_2 = \Big\{ \frac1{\sqrt{2}} \,  x_2 , x_4, \frac1{\sqrt{2}} \, x_2 \, ,
x_6 , \frac1{\sqrt{2}}\, x_8 , x_{10} , \frac1{\sqrt{2}}\, x_8 , x_{12} ,
 \frac1{\sqrt{2}}\, x_1 ,  \dots \Big\} \in \fsa \ ,
$$
with $e(\F_2 ) = 2$. In a similar way, it can be constructed frames $\F_k \in \fsa$ with
$e(\F_k ) = k$,  for every $k \in \N\cup \{ \infty\}$. Note that
$$
\F_\infty = \Big\{ \frac1{\sqrt{2}}\, x_1, x_4, \frac1{\sqrt{2}}\, x_2, x_8,
\frac1{\sqrt{2}}\, x_2 , x_{12} ,\frac1{\sqrt{2}}\, x_3, x_{16},
\frac1{\sqrt{2}}\, x_6, x_{20},\frac1{\sqrt{2}}\, x_6, \dots \Big\}  .
$$
In other words, $\F_\infty$ is the frame induced by the  bounded operator $T: \ell^2(\N) \to \hil$
given by
$$
T(e_n ) = \begin{cases} \ \  x_{4k} & \peso{if} n = 2k  \\ \ \ \frac1{\sqrt{2}}\,
x_{2k-1} & \peso{if} n = 6k-5 \\ \ \ \frac1{\sqrt{2}}\,
x_{4k-2}  & \peso{if} n = 6k-3 \\ \ \
\frac1{\sqrt{2}}\, x_{4k-2}  & \peso{if} n = 6k-1 \ .
\end{cases}
$$
Therefore $\nulli{S}{\a} = \N \cup \{0,  \infty\} \ $.  \EOE
\end{exa}

\begin{pro}\label{infi}
Let $S\in \glh^+ $ and $\c\in \ell^\infty (\N )^+$.
Assume that \fcom  \ and
$
\liminf \c < \min \ \sigma_{e}(S).
$
Then $\nulli{S}{\c}=\{\infty\}$.
\end{pro}

\bdem
 Let $\cF = \fram
 \in \fsc$, with $e(\cF ) = d$.
 By Proposition \ref{la equi} there exists  an extension $\cK = \H\oplus\H_d$ of $\H$
such that, if we denote
$$
S_1 = \begin{pmatrix} S&0\\0&0\end{pmatrix} \barr {c} \H\\\H_d \earr \in
L(\cK ) ^+ \ , \peso{then } \
\c \in \cC  \left[ \uks  \right ] \ .
$$
By Theorem  \ref{condiciones necesarias en general},
$\min \ \sigma_{e} (S_1)  \leq  \liminf \ \c \ $.
But,  if $\dim \hil_d =e(\cF ) <\infty$, then   $ \sigma_{e} (S_1)
= \sigma_{e}(S)$, which contradicts the fact that $
\liminf \c < \min \ \sigma_{e}(S)
$.
 \edem

\begin{rem} \rm  
Let $\F = \fram$ be a Parseval frame for $\hil$
(i.e. it has frame operator $S = I_\hil$). If
$
\liminf_{n \in \N } \|f_n\|  < 1 \ ,
$
then, by Proposition \ref{infi}, $e(\F ) = \infty \ $. 
This results was proved in \cite {[BCHL]}\EOE
\end{rem}

\begin{exa}

Let $\H$ be a Hilbert space with an orthonormal basis
$\cB = \{x_i\}_{i\in \N}$.
Let
$$
\a = (1, 2, 1, 2,\dots ) \ , \ \   S = M_{\cB , \a} \in \glh^+ \peso{ and } \c =
\big(\frac32 \, , \frac32 \, , \frac32 \, ,  \dots \big).
$$
 We shall show that also
$\nulli{S}{\c} = \N \cup \{0,  \infty\}$ . Note that, in this case,
$$
\alpha_-(S) = 1 <  \liminf \c = \frac 32 =  \limsup \c < 2 = \|S\|_{e} \ .
$$
Indeed, take the Riesz basis $\F_0 = \{f_n \}_{n \in \N}$ given by
$$
f_n = \begin{cases} \frac{x_n}{\sqrt{2}}+ x_{n+1}  & \peso{if $n$  is odd} \\ &\\
\frac{-x_{n-1}}{\sqrt{2}}+ x_{n}  & \peso{if $n$  is even } \end{cases} \ .
$$
It is easy to see that $\F_0 \in \fsc$. Using that
$$
\big(\ \frac32 \, , \frac32 \, , \frac32 \,  , \frac32 \ \big) \prec (2, 2, 2, 0 ) \ ,
$$
an arbitrary number of packs of four vectors with norm $\sqrt{\frac{{3}}{2}}$ associated to
packs of three even places of the diagonal of
$S$  can be interlaced into the previous construction.
Each of these packs   adds excess $1$ to the whole system.
In this way,  frames
$\F_k \in \fsc$ with $e(\F_k) = k$  can be found
for every $k \in \N \cup \{\infty\} $.
\EOE
\end{exa}

\begin{rem}\rm
Let $S \in  \matinv^+$ and $\c \in \ell^\infty (\M)^+ $.
If \fcom , then $\nulli{S}{\c} = \{ |\M| - n \}$.
Nevertheless, if $k>n$, $\c = (1, \dots , 1) \in
\C^k$ and $S = \frac kn I \in \mat$, then $\fsc$ is the set of spherical tight frames
of $k$ elements in $\cene$. Dykema, Freeman, Korleson, Larson, Ordower
and Weber \cite{Dyk} have shown that, in this case, $\fsc$
has a rich geometrical structure, with several orbits of qualitatively
different elements. \EOE
\end{rem}

\subsection*{Acknowledgments:}
We would like to acknowledge Professor
G. Corach, who shared with us fruitful discussions concerning
the problems included in this article. We also thank the referee 
for calling our attention to the papers of Kadison \cite{Kad1}, \cite{Kad2}.

\end{document}